\documentclass[a4paper,11pt]{article}

\input{diagrams}

\newcommand{\op}{\ensuremath{\mbox{\hspace{1pt}{\scriptsize op}}}}
\newcommand{\set}{\ensuremath{\mbox{\bfseries Set}}}
\newcommand{\Set}{\ensuremath{\mbox{\bfseries Set}}}
\newcommand{\cat}[1]{\ensuremath{\mbox{\bfseries {\upshape {#1}}}}}
\newcommand{\free}[1]{\ensuremath{{\mathcal F} \hspace{-1.5pt}
    {#1}^{\op}}}

\newcommand{\elt}[1]{\ensuremath{\mbox{\upshape elt} \hspace{1pt} {#1}}}
\newcommand{\sfree}[1]{\ensuremath{{\mathcal F} \hspace{-1.5pt}
    {#1}^{\mbox{\hspace{1pt}{\tiny op}}}}}
\newcommand{\cl}[1]{\ensuremath{\mathcal {#1}}}
\newcommand{\bb}[1]{\ensuremath{\mathbb {#1}}}

\newcommand{\xxx}{\ensuremath{
    \mbox{elt}\hspace{1pt}Q \times_{{\mathcal
    F}\hspace{-1pt}S^{\mbox{{\tiny op}}}}{\mathcal
    F}\hspace{-1pt}X^{\mbox{{\scriptsize op}}}}}

\newcommand{\xx}[3]{\ensuremath{
\mbox{elt}\hspace{1pt}{#1} \times_{{\mathcal
F}\hspace{-1pt}{#2}^{\mbox{{\tiny op}}}}{\mathcal
F}\hspace{-1pt}{#3}^{\mbox{{\scriptsize op}}}}}

\newcommand{\xv}[1]{\ensuremath{{#1}_1,\ldots,{#1}_n}}

\newcommand{\vect}[3]{\renewcommand{\arraystretch}{0.5}\ensuremath{
    \left(\begin{array}{@{} c @{\hspace{0pt}}l @{}} {#1} \\
    \downarrow & {#2} \\ {#3} \end{array}\right)}}

\newcommand{\tomspan}[8]
{\begin{picture}(80,50)      %
\put(8,10){${#2}\vect{#3}{#4}{#1}$}  %  bottom left
\put(29.5,35){$\vect{#5}{#6}{#1}$}  %  top
\put(47,10){$\vect{#3}{#4}{#1}$}  %  bottom right

\put(29,35){\vector(-2,-3){10}}  %
\put(40,35){\vector(2,-3){10}}  %
\put(24,30){\makebox(0,0)[r]{${#7}$}} %  left
\put(45,30){\makebox(0,0)[l]{${#8}$}}%  right

\end{picture}}

\newcommand{\lra}{\ensuremath{\longrightarrow}}
\newcommand{\map}[1]{\ensuremath{\stackrel{{#1}}{\lra}}}

\newcommand{\sm}{symmetric multicategory{}}
\newcommand{\sms}{symmetric multicategories{}}
\newcommand{\topd}{$(\cl{E},T)$-operad}
\newcommand{\tmcat}{$(\cl{E},T)$-multicategory}
\newcommand{\tmcats}{$(\cl{E},T)$-multicategories{}}
\newcommand{\tqopd}{$(\cl{E}_Q,T_Q)$-operad}

\newcommand{\st}{\ensuremath{(\cl{E},T)}}
\newcommand{\mcat}{multicategory }

\newcommand{\sqtq}{\ensuremath{(\cl{E}_Q, T_Q){}}}
\newcommand{\monic}{\ensuremath{\succ \hspace{-3pt}
    \longrightarrow}}
\newcommand{\tqplus}{\ensuremath{(\cl{E}_{Q^+},T_{Q^+})}}
\newcommand{\tqprime}{\ensuremath{({\cl{E}_Q}', {T_Q}')}}
\newcommand{\tpr}{\ensuremath{{T_Q}'}}
\newcommand{\tpl}{\ensuremath{T_{Q^+}}}
\newcommand{\cspl}{\ensuremath{\cl{E}_{Q^+}}}
\newcommand{\cspr}{\ensuremath{{\cl{E}_Q}'}}

\newcommand{\iso}{isomorphism{}}

\newcommand{\tck}{\ensuremath{T(C^{(k)})}}
\newcommand{\catck}{\ensuremath{\bb{C}^{(k)}}}
\newcommand{\ck}{\ensuremath{C^{(k)}}}
\newcommand{\id}{\ensuremath{\mbox{\em id}}}

\newtheorem{theorem}{Theorem}[section]
\newtheorem{proposition}[theorem]{Proposition}
\newtheorem{corollary}[theorem]{Corollary}

\newtheorem{definition}[theorem]{Definition}

\newenvironment{prf}{\vspace{2ex}\begin{sloppypar}{\noindent\upshape
{\bfseries Proof. }}} {{\hspace*{\fill}
$\Box$}\end{sloppypar}\vspace{2ex}}

\newcommand{\numarabic}{\renewcommand{\labelenumi}{\arabic{enumi})}}

\newcommand{\pica}{\begin{center} \input}
\newcommand{\picz}{\end{center}}
\newcommand{\length}[1]{\setlength{\unitlength}{#1}}

\newlength{\leng}
\newlength{\fontleng}
\newcommand{\sunit}{\setlength{\unitlength}{1mm}}

\sunit

\sunit

\sunit

\sunit

\sunit

\sunit

\sunit

\sunit

\sunit

\usepackage{amsfonts}
\usepackage{amssymb}
\usepackage{amsmath}
\usepackage{fancyheadings}
\usepackage{latexcad}

\begin{document}

\title{Weak $n$-categories: comparing opetopic foundations}
\author{Eugenia Cheng\\ \\Department of Pure Mathematics, University
of Cambridge\\E-mail: e.cheng@dpmms.cam.ac.uk}
\date{October 2002}
\maketitle

\begin{abstract}
We define the category of tidy symmetric multicategories.  We
construct for each tidy symmetric multicategory $Q$ a cartesian
monad $(\cl{E}_Q,T_Q)$ and extend this assignation to a functor.  We exhibit a relationship between the
slice construction on symmetric multicategories, and the `free
operad' monad construction on suitable monads.  We use this to
give an explicit description of the relationship between
Baez-Dolan and Leinster opetopes.
\end{abstract}

\setcounter{tocdepth}{3}
\tableofcontents

\section*{Introduction}
\addcontentsline{toc}{section}{Introduction}

The problem of defining a weak $n$-category has been approached in
various different ways (\cite{bd1}, \cite{hmp1}, \cite{lei1},
\cite{pen1}, \cite{bat1}, \cite{tam1}, \cite{str2}, \cite{may1},
\cite{lei7}), but so far the relationship between
these approaches has not been fully understood.  The subject of
the present paper is the relationship between the approaches given
in \cite{bd1} and \cite{lei1}.  This work is a continuation of the
work of \cite{che7}, in which we exhibited a relationship between
the approaches of \cite{bd1} and \cite{hmp1}.  The effect of the
present paper is therefore to establish a relationship between all
three approaches.  We refer the reader to Section~1 of \cite{che7} for an
overview of the subject.

In each of the cases in question the definition has two components. First, the
language for describing $k$-cells is set up.  Then, a concept of universality
is introduced, to deal with composition and coherence.  Any comparison of these
approaches must therefore begin at the construction of $k$-cells, and in this
paper, as in \cite{che7}, we restrict our attention to this process.  This, in
the terminology of Baez and Dolan, is the theory of opetopes.

The idea behind the opetopic approaches is to construct opetopes
as the underlying shapes of cells.  These are described using the
language of multicategories.  $k$-cells are constructed from
formal composites or `pasting diagrams' of
($k$\hspace{1pt}-1)-cells, so a $k$-cell is considered as an arrow
from a list of its constituent ($k$\hspace{1pt}-1)-cells to a
($k$\hspace{1pt}-1)-cell. This raises an immediate question: in
what order should we list the constituent
($k$\hspace{1pt}-1)-cells?  As Leinster points out in \cite{lei1},
there is in general no ordering that is stable under composition.
The different approaches we consider here arise from, essentially,
different ways of dealing with this problem.

John Baez and James Dolan tackle the problem by including all
possible orderings, giving rise to a symmetric action. In
\cite{bd1}, they give a definition of weak $n$-categories based on
operads (symmetric multicategories), opetopes and opetopic sets.

In \cite{hmp1}, Claudio Hermida, Michael Makkai and John Power
begin an explicitly analogous definition, based on (generalised)
multicategories, multitopes and multitopic sets.  This version
arises from tackling the above problem by choosing one ordering
for each cell; the resulting complications are dealt with by
generalising the notion of composition.

In \cite{che7} the relationship between
these approaches is made explicit, as far as the notions of
opetopes and multitopes.

First we exhibit an embedding
	\[\xi:\mbox{{\bfseries GenMulticat}} \hookrightarrow
\mbox{{\bfseries SymMulticat}}\]
of the category of generalised multicategories into the category of
symmetric multicategories.  We then examine the slicing process, used
for constructing $k$-cells from $(k-1)$-cells, and we show that the functor
$\xi$ ``commutes'' with slicing (up to equivalence).  Finally, since
opetopes and multitopes arise from iterated slicing, we deduce that
the respective categories of $k$-opetopes and $k$-multitopes are equivalent.

In fact, we do not use the definition of opetopes precisely as
given in \cite{bd1}, but rather, we develop a generalisation of
the notion along lines which Baez and Dolan began but chose to
abandon, for reasons unknown to the present author.  

Motivated by this work, in the present paper we continue to use
the generalisation of the opetopic approach rather than the
approach precisely as given in \cite{bd1}.

In \cite{lei1}, Tom Leinster gives an approach to the construction
which is quite different; he tackles the problem of ordering
constituent cells by not ordering them at all.  Instead of attempting
to squash the constituent cells into a line for the purpose of
listing them as the source of an arrow, they are allowed to remain in
their natural `positions', and the notion of {\em arrow} is
generalised instead.  So this approach is based on $({\mathcal
E},T)$-multicategories; these structures were defined by Burroni
(\cite{bur1}) and have also been treated by Hermida (\cite{her1}).
Here, the source of an arrow is not necessarily given as a list of
objects, but as a structure given by some cartesian monad $T$.

The role that these (even more generalised) multicategories plays
is not explicitly analogous to that of operads and multicategories
in the opetopic and multitopic versions respectively, so the
comparison is more subtle.  In fact, rather than comparing the
role of symmetric multicategories with that of $({\mathcal
E},T)$-multicategories, we compare it with the role of cartesian
monads.  This is the subject of Section~{\ref{taczeta}}.

To study this relationship, we can restrict our attention to a certain
kind of particularly well behaved symmetric multicategory.  
We call a symmetric multicategory {\em tidy} if it is freely
symmetric and has a category of objects equivalent to a discrete
category.  We observe (\cite{che7}) that the multicategories
involved in the construction of Baez-Dolan opetopes are all tidy;
also, a \sm\ is equivalent to one in the image of $\xi$ if and
only if it is tidy. We write \cat{TidySymMulticat} for the full
subcategory of \cat{SymMulticat} whose objects are tidy \sms.

In Section \ref{taczeta} we construct a functor
    \[\zeta : \cat{TidySymMulticat} \lra \cat{CartMonad}\]
where \cat{CartMonad} is the category of cartesian monads and
cartesian monad opfunctors.  Roughly, given a tidy \sm\ $Q$ we
construct a cartesian monad \sqtq\ which acts on sets of `labelled
$Q$-objects' to give sets of `source-labelled $Q$-arrows'.

The idea is that, for an \tmcat, much information about an arrow
is given by its domain, that is, the action of $T$; in the
Baez-Dolan setting the domain of an arrow is just a list of
objects, and the information is captured elsewhere. So, the
functor part of $T_Q$ is constructed from the collection of arrows
itself, the unit from the identities, and multiplication from the
reduction laws of $Q$.

In Section \ref{opes} we examine the construction of opetopes.  First
we examine the process of constructing
$k$-cells from ($k$\hspace{1pt}-1)-cells.  In \cite{bd1} Baez and
Dolan define the `slicing' process for this purpose.  Leinster
does define a slicing process on \tmcats, but since we are
considering a comparison between symmetric multicategories and
{\em monads}, we seek an analogous process defined on these
monads, rather than on \tmcats.  This process is the `free \topd\
monad' construction, defined on {\em suitable} monads
(\cite{lei4}).

Given a suitable monad \st, the monad $\st' = (\cl{E}',T')$ is
defined to be the free \topd\ monad.  We show that
    \[\zeta(Q^+) \cong {\zeta(Q)}'.\]
In this sense, the processes are analogous.

Finally, we apply these results to the
construction of opetopes. Having established a relationship
between the underlying theories, it is straightforward to compare
these constructions. In the Baez-Dolan setting, the category of
$k$-opetopes is defined to be the object-category of $I^{k+}$,
where $I$ is the symmetric multicategory with only one object and
one (identity) arrow.

Leinster gives a construction of `opetopes' with a role analogous to
that of Baez-Dolan opetopes, based on a series $(\set/S_n,T_n)$ of
cartesian monads.  We show that, for each $n \geq 0$ \[\zeta(I^{n+})
\cong (\set/S_n,T_n)\] and deduce that \[o(I^{n+}) \simeq S_n\] for
each $n \geq 0$, where $o$ denotes the object-category.  Informally,
we see that Baez-Dolan opetopes and Leinster opetopes are the same up
to isomorphism.

Throughout this paper we repeatedly find that the details of
proofs are fiddly but uninteresting.  So we include some informal
comments about how the constructions may be interpreted, as a
gesture towards demonstrating that the notions are in fact
naturally arising.  The aim is to shed some light on the
relationship between the various structures involved, a
relationship which has previously remained unclear.

Finally, we note that some of the constructions given in
\cite{lei1} and \cite{lei4} may be treated in slightly greater
detail in \cite{lei5}, in particular the free multicategory and
slicing constructions.

\subsubsection*{Terminology}

\renewcommand{\labelenumi}{\roman{enumi})}
\begin{enumerate}

\item Since we are concerned chiefly with {\em weak}
$n$-categories, we follow Baez and Dolan (\cite{bd1}) and omit the
word `weak' unless emphasis is required; we refer to {\em strict}
$n$-categories as `strict $n$-categories'.

\item We use the term `weak $n$-functor' for an $n$-functor where
functoriality holds up to coherent isomorphisms, and `lax' functor
when the constraints are not necessarily invertible.

\item In \cite{bd1} Baez and Dolan use the terms `operad' and
`types' where we use `multicategory' and `objects'; the latter
terminology is more consistent with Leinster's use of `operad' to
describe a multicategory whose `objects-object' is 1.

\item In \cite{hmp1} Hermida, Makkai and Power use the term
`multitope' for the objects constructed in analogy with the
`opetopes' of \cite{bd1}.  This is intended to reflect the fact
that opetopes are constructed using operads but multitopes using
multicategories, a distinction that we have removed by using the
term `multicategory' in both cases.  However, we continue to use
the term `opetope' and furthermore, use it in general to refer to
the analogous objects constructed in each of the three theories.
Note also that Leinster uses the term `opetope' to describe
objects which are analogous but not {\em a priori} the same; we
refer to these as `Leinster opetopes' if clarification is needed.

\item We follow Leinster and use the term `$(\cl{E},T)$-multicategory' for the
notion defined by Burroni (\cite{bur1}) as `$T$-category' (in
French).

\item We regard sets as sets or discrete categories with no
notational distinction.

\end{enumerate}

\bigskip {\bfseries Acknowledgements}

This work was supported by a PhD grant from EPSRC.  I would like to thank
Martin Hyland and Tom Leinster for their support and guidance.

\section{The theory of multicategories}
\label{newopes}

Opetopes are described using the language of multicategories.  In
each of the three theories of opetopes in question, a different
underlying theory of multicategories is used.  In this section we
examine the underlying theories, and we construct a way of
relating these theories to one another; this relationship provides
subsequent equivalences between the definitions. We adopt a concrete
approach here; certain aspects of the definitions suggest a more
abstract approach but this will require further work beyond the scope
of this work.

%In this section we compare the three theories of multicategories
%used in the definitions of \cite{bd1}, \cite{hmp1} and \cite{lei1}
%respectively, and show that they are equivalent.

\label{tacmcats}

\subsection{Symmetric multicategories} \label{tacsym}
\numarabic

By `symmetric multicategory' we mean symmetric multicategory with a
category of objects.  Such structures are defined (as `operads') in
\cite{bd1}; see also \cite{che7}.  We recall the definition here.

We write $\cl{F}$ for the `free symmetric strict monoidal
category' monad on \cat{Cat}, and ${\mathbf S}_k$ for the group of
permutations on $k$ objects; we also write $\iota$ for the
identity permutation.

\begin{definition} A {\em symmetric multicategory} $Q$ is given by
the following data
\begin{enumerate}

\item A category $o(Q)=\bb{C}$ of objects. We refer to \bb{C} as the
{\em object-category}, the morphisms of \/ ${\mathbb C}$ as {\em
object-morphisms}, and if \/ ${\mathbb C}$ is discrete, we say
that $Q$ is {\em object-discrete}.

\item For each $p\in \cl{F}{\mathbb C}^{\mbox{\scriptsize op}} \times
{\mathbb C}$, a set $Q(p)$ of arrows. Writing
    \[p=(x_1, \ldots ,x_k;x),\]
an element $f \in Q(p)$ is considered as an arrow with source and
target given by
\begin{eqnarray*}
    s(f) &=& (x_1,\ldots ,x_k)\\
    t(f) &=& x
\end{eqnarray*}
and we say $f$ has {\em arity} $k$.  We may also write $a(Q)$ for
the set of all arrows of $Q$.

\item For each object-morphism $f:x \longrightarrow y$,
an arrow $\iota(f) \in Q(x;y)$.  In particular we write
$1_x=\iota(1_x)\in Q(x;x)$.

\item Composition: for any $f \in Q(x_1, \ldots ,x_k;x)$ and
$g_i \in Q(x_{i1}, \ldots ,x_{im_i};x_i)$ for $1 \le i \le k$, a
composite
\[f \circ (g_1, \ldots ,g_k) \in Q(x_{11}, \ldots ,x_{1m_1},
\ldots ,x_{k1}, \ldots ,x_{km_k};x)\]

\item Symmetric action: for each permutation $\sigma \in {\mathbf S}_k$, a
    map
    \[\begin{array}{rccc} \sigma : & Q(x_{1}, \ldots ,x_{k};x)
    & \lra & Q(x_{\sigma (1)}, \ldots ,x_{\sigma (k)};x) \\
    & f & \longmapsto & f \sigma
    \end{array}\]

\end{enumerate}

\noindent satisfying the following axioms:

\begin{enumerate}

\item Unit laws: for any $f \in Q(x_1, \ldots ,x_m;x)$, we have
\[1_x \circ f = f = f \circ (1_{x_1}, \ldots, 1_{x_m})\]

\item Associativity: whenever both sides are defined,
 \[ \begin{array}{c} f \circ(g_1 \circ (h_{11}, \ldots , h_{1m_1}),
    \ldots , g_k \circ (h_{k1}, \ldots , h_{km_k})) = \mbox{\hspace{5em}}\\
    \mbox{\hspace*{5em}} (f \circ (g_1, \ldots , g_k)) \circ (h_{11},
    \dots ,h_{1m_1}, \ldots ,h_{k1}, \ldots ,h_{km_k}) \end{array}\]

\item For any $f \in Q(x_1, \ldots ,x_m;x)$ and
$\sigma , \sigma' \in {\mathbf S}_k$,
    \[(f \sigma)\sigma' = f(\sigma\sigma')\]

\item For any $f \in Q(x_1, \ldots ,x_k;x)$, $g_i \in Q(x_{i1},
\ldots ,x_{im_i};x_i)$ for $1 \le i \le k$, and $\sigma \in
{\mathbf S}_k$, we have
    \[(f \sigma) \circ (g_{\sigma (1)}, \ldots
    ,g_{\sigma(k)}) =
    f \circ (g_1, \ldots , g_k) \cdot \rho(\sigma)\]
where $\rho :{\mathbf S}_k \longrightarrow {\mathbf S}_{m_1 +
\ldots + m_k}$ is the obvious homomorphism.

\item For any $f \in Q(x_1, \ldots ,x_k;x)$, $g_i \in Q(x_{i1},
\ldots ,x_{im_i};x_i)$, and $\sigma_i \in {\mathbf S}_{m_i}$ for
$1 \le i \le k$, we have
    \[f \circ (g_1 \sigma_1, \ldots, g_k \sigma_k) =
    (f \circ(g_1, \ldots, g_k))\sigma\]
where $\sigma \in {\mathbf S}_{m_1 + \dots + m_k}$ is the
permutation obtained by juxtaposing the $\sigma_i$.

\item $\iota(f\circ g) = \iota(f) \circ \iota(g)$

\end{enumerate}
\end{definition}

%\label{tacspan} A symmetric multicategory $Q$ may be thought of as
%a functor
%   \[Q:\cl{F}\bb{C}^{\op} \times {\mathbb C}
%    \longrightarrow \cat{Set}\]
%with some extra structure.

%In a more abstract view, we would expect \cl{F} to be a 2-monad on
%the 2-category \cat{Cat}, which lifts via a generalised form of
%distributivity to a bimonad on \cat{Prof}, the bicategory of
%profunctors.  Then the Kleisli bicategory for this bimonad should
%have as objects small categories, and its 1-cells should be
%essentially profunctors of the form $\cl{F}\bb{C}
%\makebox[0pt][l]{\hspace{1em}$\mid$}\longrightarrow \bb{D}$ in the
%opposite category. However, the calculations involved in this
%description are intricate and require further work.

%In this abstract view, a symmetric multicategory $Q$ would then be
%a monad in this bicategory.  Arrows and symmetric action (Data 2,
%5) are given by the action of $Q$, identities (Data 3) by the unit
%of the monad and composition (Data 4) by the multiplication for
%the monad.

%1

\begin{definition} Let $Q$ and $R$ be symmetric multicategories
with object-categories \bb{C} and \bb{D} respectively. A {\em
morphism of symmetric multicategories} $F:Q \longrightarrow R$ is
given by

\begin{itemize}

\item A functor $F=F_0:{\mathbb C} \longrightarrow {\mathbb D}$

\item For each arrow $f \in Q(x_1, \ldots ,x_k;x)$ an
arrow $Ff \in R(Fx_1, \ldots ,Fx_k; Fx)$

\end{itemize}

\noindent satisfying

\begin{itemize}

\item $F$ preserves identities: $F(\iota(f)) = \iota(Ff)$ so in
particular $F(1_x) = 1_{Fx}$

\item $F$ preserves composition: whenever it is defined
\[F(f \circ (g_1,\ldots,g_k)) = (Ff \circ (Fg_1, \ldots ,
Fg_k))\]

\item $F$ preserves symmetric action: for each $f \in Q(x_1,
\ldots , x_k;x)$ and $\sigma \in {\mathbf S}_k$
\[F(f\sigma) = (Ff) \sigma \]

\end{itemize}
\end{definition}

\noindent Composition of such morphisms is defined in the obvious
way, and there is an obvious identity morphism $1_Q:Q
\longrightarrow Q$. Thus symmetric multicategories and their
morphisms form a category {\bf SymMulticat}.

\begin{definition} A morphism $F:Q \longrightarrow R$ is an {\em equivalence} if and
only if the functor $F_0:{\mathbb C} \longrightarrow {\mathbb D}$
is an equivalence, and $F$ is full and faithful. That is, given
objects $x_1, \ldots, x_m, x$ the induced function
    \[F: Q(x_1, \ldots, x_m; x) \lra R(Fx_1, \ldots, Fx_m; Fx)\]
is an isomorphism.
\end{definition}

%\noindent The following observations will be useful later.  If the
%functor $F_0$ is an equivalence and $R$ is object-discrete, then
%\begin{enumerate}
%
%\item $F_0$ must send all object-morphisms to identities in
%${\mathbb D}$.  So for any object-morphism $f:x \longrightarrow
%y$, we must have $F_0(x)=F_0(y)$.
%
%\item $F_0$ is full, so if $F_0(x)=F_0(y)$ then there must exist
%object-morphisms
%\[\begin{array}{c}
%    f:x \longrightarrow y\\
%    g:y \longrightarrow x. \end{array}\]
%
%\item $F_0$ is faithful, so if $F_0(x)=F_0(y)$ then there is at
%most one object-morphism $x \longrightarrow y$.
%
%\end{enumerate}
%
%\noindent So in these circumstances all object-morphisms must be
%unique isomorphisms, and
%    \[F_0(x)=F_0(y) \iff x \cong y \mbox{ in } {\mathbb C}.\]
%

\subsection{\tmcats}
\label{tactmcat}

In \cite{lei1} opetopes are constructed using
$(\cl{E},T)$-multicategories.  These are defined by Burroni in
\cite{bur1} as `$T$-categories'.

\begin{definition} Let $T$ be a cartesian monad on a cartesian category
\cl{E}.  An {\em $(\cl{E},T)$-multicategory} is given by an
`objects-object' $C_0$ and an `arrows-object' $C_1$, with a
diagram
    \[TC_0 \stackrel{d}{\longleftarrow} C_1 \map{c} C_0\]
in \cl{E} together with maps $C_0 \map{\mbox{{\tiny\em ids}}} C_1$
and $C_1 \circ C_1 \map{\mbox{{\tiny\em comp}}} C_1$ satisfying
associative and identity laws.  (See \cite{lei5} for full
details.)
\end{definition}

We write {\bfseries CartMonad} for the category of cartesian
monads and cartesian monad opfunctors.  A cartesian monad
opfunctor
    \[(U,\phi):({\mathcal E}_1,T_1) \longrightarrow ({\mathcal
    E}_2,T_2)\]
consists of \begin{itemize}
\item a functor $U: \cl{E}_1 \longrightarrow \cl{E}_2$
preserving pullbacks
\item a cartesian natural transformation $\phi: UT_1
\longrightarrow T_2U$, that is, a natural transformation whose
naturality squares are pullbacks
\end{itemize} satisfying certain axioms (see \cite{str1} and \cite{lei4}
for full definitions).

\subsection{Relationship between symmetric multicategories and
cartesian monads} \label{taczeta}

The respective roles of multicategories in the Baez-Dolan and
Leinster approaches are not explicitly analogous.  In this section
we  exhibit instead a correspondence between certain symmetric
multicategories and certain cartesian monads, by constructing a
functor
    \[\zeta : \mbox{\cat{TidySymMulticat}} \longrightarrow
    \mbox{\cat{CartMonad}}.\]
This is enough since we will see that all the symmetric
multicategories involved in the construction of opetopes are tidy.

We begin by defining the action of $\zeta$ on objects; so for any
tidy symmetric multicategory $Q$, we construct a cartesian monad
$\zeta(Q) = (\cl{E}_Q, T_Q)$, say. Informally, the idea behind
this construction is that $T_Q$ should encapsulate information
about the arrows of $Q$.  The functor part is constructed to give
the arrows themselves, the unit to give the identities, and
multiplication the reduction laws (composites).

Write $o(Q) = \bb{C}$.  $Q$ is tidy, so $\bb{C} \simeq S$, say,
where $S$ is a discrete category. For various of the constructions
which follow, we assume that we have chosen a specific functor $S
\map{\sim} \bb{C}$. However, when isomorphism classes are taken
subsequently, we observe that the construction in question does
not depend on the choice of this functor.

Put $\cl{E}_Q = \cat{Set}/S$ and observe immediately that this is
cartesian.  (This is sufficient here, though of course \set/S has
much more structure than this.)

Informally, an element $(X,f)=(X \map{f} S)$ of $\set/S$ may be
thought of as a system for labelling $Q$-objects with `compatible'
elements of X; each `label' is compatible with an isomorphism
class of $Q$-objects. Then the action of $T_Q$ assigns compatible
labels to the source elements of $Q$-arrows in every way possible;
the target is not affected. The resulting set of `source-labelled
$Q$-arrows' is itself made into a set of labels by regarding each
arrow as a `label' for its target.

We now give the formal definition of the functor $T_Q: \cl{E}_Q
\lra \cl{E}_Q$. For the action on object-categories, consider
$(X,f) = (X \stackrel{f}{\longrightarrow} S) \in \set/S$. We have
the following composite functor
    \[\elt{Q} \stackrel {s}
    {\longrightarrow}\cl{F}\bb{C}^{\op} \stackrel
    {\sim}{\longrightarrow} \free S\]
where \cl{F} denotes the free symmetric strict monoidal category
monad on \cat{Cat}, and $s$ and $t$ the source and target
functions respectively. Consider the pullback

\begin{center}
\setlength{\unitlength}{0.8mm}
\begin{picture}(50,35)(5,5)      %

\put(10,10){\makebox(0,0){\elt{Q}}}  %  bottom left
\put(10,35){\makebox(0,0){$\cdot$}}  %  top left
\put(47,35){\makebox(0,0){\free{X}}}  % top right
\put(47,10){\makebox(0,0){\free{S} .}}  % bottom right

\put(13,35){\vector(1,0){25}}  %     top
\put(16,10){\vector(1,0){23}}  %     bottom
\put(10,32){\vector(0,-1){17}} %     left
\put(45,32){\vector(0,-1){17}} %     right

\put(8,23){\makebox(0,0)[r]{}} %           left
\put(47,23){\makebox(0,0)[l]{\free{f}}}%   right
\put(27,37){\makebox(0,0)[b]{}}%           top
\put(27,8){\makebox(0,0)[t]{}} %           bottom

\end{picture}
\end{center} Since $Q$ is tidy, \elt{Q} is equivalent to a discrete
category, and so too is the above pullback; so we have \[\elt{Q}
\times_{\sfree{S}} \free{X} \simeq X',\] say, where $X'$ is
discrete. Put $T_Q (X,f) = (X',f')$ where $f'$ is the composite
    \[X' \stackrel{\sim}{\longrightarrow} \elt{Q}
    \times_{\sfree{S}} \free{X} \longrightarrow \elt{Q} \stackrel
    {t}{\longrightarrow} {\mathbb C} \stackrel
    {\sim}{\longrightarrow} S.\]
This is well-defined since if $(\alpha,\underline{x})
\cong(\alpha',\underline{x}') \in \elt{Q} \times_{\sfree{S}}
\free{X}\ $ then certainly $\alpha \cong \alpha' \in \elt{Q}$ and
so $t(\alpha) \cong t(\alpha') \in {\mathbb C}$.

%7

We now define the action of $T_Q$ on morphisms.  A morphism

\begin{center}
\setlength{\unitlength}{1mm}
\begin{picture}(37,22)(2,12)      %

\put(20,12){\makebox(0,0){$S$}}  %  bottom
\put(5,30){\makebox(0,0){$X$}}  %  top left
\put(35,30){\makebox(0,0){$Y$}}  %  top right

\put(8,30){\vector(1,0){24}}  %  top
\put(5,27){\vector(1,-1){13}}  %  left
\put(34,27){\vector(-1,-1){13}} %  right

\put(30,20){\makebox(0,0)[l]{$g$}} %  right
\put(9,20){\makebox(0,0)[r]{$f$}}%  left
\put(20,31){\makebox(0,0)[b]{$h$}}%  top

\end{picture}
\end{center}in $\set/S$ induces a functor
    \[\elt{Q} \times_{\sfree{S}} \free{X}
    \longrightarrow \elt{Q} \times_{\sfree{S}} \free{Y}\]
which, by construction, makes the following diagram commute:

\begin{center}
\setlength{\unitlength}{1mm}
\begin{picture}(42,25)(0,10)      %

\put(20,11){\makebox(0,0){\elt{Q}}}  %  bottom
\put(0,30){\makebox(0,0){$\elt{Q} \times_{\sfree{S}} \free{X}$}}
\put(40,30){\makebox(0,0){$\elt{Q} \times_{\sfree{S}} \free{Y}$}} %top right

\put(15,30){\vector(1,0){10}}  %  top
\put(5,27){\vector(1,-1){13}}  %  left
\put(34,27){\vector(-1,-1){13}} %  right

\put(30,20){\makebox(0,0)[l]{}} %  right
\put(9,20){\makebox(0,0)[r]{}}%  left
\put(20,31){\makebox(0,0)[b]{}}%  top

\end{picture}
\end{center}giving a morphism

\begin{center}
\setlength{\unitlength}{1mm}
\begin{picture}(36,25)(3,10)       %

\put(20,12){\makebox(0,0){$S$}}  %  bottom
\put(5,30){\makebox(0,0){$X'$}}  %  top left
\put(35,30){\makebox(0,0){$Y'$}}  %  top right

\put(8,30){\vector(1,0){24}}  %  top
\put(5,27){\vector(1,-1){13}}  %  left
\put(34,27){\vector(-1,-1){13}} %  right

\put(30,20){\makebox(0,0)[l]{$g'$}} %  right
\put(9,20){\makebox(0,0)[r]{$f'$}}%  left
\put(20,31){\makebox(0,0)[b]{$h'$}}%  top

\end{picture}
\end{center}in $\Set/S$.  We define $T_Q$ on morphisms by $T_Q(h) = h'$.
$T_Q$ is clearly functorial; we now show that it inherits a
cartesian  monad structure from the identities and composition of
$Q$.  For convenience we write $\cl{E}_Q=\cl{E}$ and $T_Q=T$.

\begin{itemize}
\item unit
\end{itemize}

We seek a natural transformation $\eta: 1_{\cl{E}} \Longrightarrow
T$, so with the above notation we need components
    \[\eta_{(X,f)}:(X,f) \lra (X',f').\]

Given $(X,f) \in \set/S$, we have a functor $X \longrightarrow
\elt{Q}$ \ given by the composite
    \[X \stackrel{f}{\longrightarrow} S
    \stackrel{\sim}{\longrightarrow} \bb{C}
    \stackrel{1_{\_}}{\longrightarrow} \elt{Q}.\]
We also have a functor $X \longrightarrow \free{X}$ given by the
unit of the monad ${\mathcal F}$. These induce a functor
    \[X \longrightarrow \xxx \]
and we define the component $\eta_{(X,f)}$ to be the composite
    \[X \longrightarrow \xxx \stackrel{\sim}{\longrightarrow}X'.\]
Explicitly, $\eta_{(X,f)}$ acts as follows.  We have
$\eta_{(X,f)}(x) = [(1_c,x)]$, the isomorphism class of
    \[(1_c,x) \in \elt{Q} \times_{\cl{F}S^{\op}} \cl{F}X^{\op}.\]
So $(1_c, x)$ is an ``identity labelled by $x$'', where $c \in
\bb{C}$ is any object in the isomorphism class $fx$.  We can see
explicitly that this is well defined since if $c \cong c'$ we have
$1_c \cong 1_{c'} \in \elt{Q}$ and thus
    \[[(1_c, x)] = [(1_{c'}, x)].\]

The following diagram commutes

\setlength{\unitlength}{1mm}

\begin{center}
\begin{picture}(40,25)(0,10)      %

\put(20,12){\makebox(0,0){$S$}}  %  bottom
\put(5,30){\makebox(0,0){$X$}}  %  top left
\put(35,30){\makebox(0,0){$X'$}}  %  top right

\put(8,30){\vector(1,0){24}}  %  top
\put(5,27){\vector(1,-1){13}}  %  left
\put(34,27){\vector(-1,-1){13}} %  right

\put(30,20){\makebox(0,0)[l]{$f'$}} %  right
\put(9,20){\makebox(0,0)[r]{$f$}}%  left
\put(20,31){\makebox(0,0)[b]{$\eta_{(X,f)}$}}%  top

\end{picture}
\end{center}

\noindent so $\eta_{(X,f)}$ is indeed a morphism $(X,f)
\longrightarrow T(X,f) \in \set/S$.

%8

Next we show that the components $\eta_{(X,f)}$ satisfy
naturality; so we show that for any morphism $h:(X,f)
\longrightarrow (Y,g) \in \set/S$ the following diagram commutes

\setlength{\unitlength}{1mm}

\begin{center}
\begin{picture}(40,35)(7,5)      %

\put(10,10){\makebox(0,0){$Y$}}  %  bottom left
\put(10,35){\makebox(0,0){$X$}}  %  top left
\put(45,35){\makebox(0,0){$X'$}}  %  top right
\put(45,10){\makebox(0,0){$Y'$}}  %  bottom right

\put(13,35){\vector(1,0){29}}  %  top
\put(13,10){\vector(1,0){29}}  %  bottom
\put(10,32){\vector(0,-1){19}} %  left
\put(45,32){\vector(0,-1){19}} %  right

\put(8,23){\makebox(0,0)[r]{$h$}} %  left
\put(47,23){\makebox(0,0)[l]{$h'$}}%  right
\put(27,37){\makebox(0,0)[b]{$\eta_{(X,f)}$}}%  top
\put(27,8){\makebox(0,0)[t]{$\eta_{(Y,g)}$}} %  bottom

\end{picture}
\end{center}

\noindent This follows from the construction of $\eta$, and
naturality of the unit for \cl{F}; alternatively, we see that on
elements, the right-ish leg gives
    \[x \longmapsto [(1_c,x)] \longmapsto [(1_c,hx)]\]
with $c$ in the isomorphism class $fx$, and the left-ish leg
gives
     \[x \longmapsto hx \longmapsto [(1_{c'},hx)]\]
with $c'$ in the isomorphism class $ghx$.  But $gh=f$ since
$h:(X,f) \longrightarrow (Y,g)$, so $c' \cong c$ and
$[(1_{c'},hx)] = [(1_c,hx)]$.

It also follows from the construction of $\eta$ that the square is
a pullback; it is similarly easily seen by considering elements.

\begin{itemize}
\item multiplication
\end{itemize}

%fff

We seek a natural transformation $\mu : T^2 \Longrightarrow T$.
Consider $(X,f) \in \set/S$.  Then by definition
\[\begin{array}{rl}
    & X' \simeq \xxx = A, \mbox{\ say}\\
    \mbox{and} & X'' \simeq \xx{Q}{S}{X'} = B, \mbox{\ say}.
\end{array}\]
We construct a commutative square

\begin{center}
\setlength{\unitlength}{0.8mm}
\begin{picture}(50,35)(5,5)      %

\put(10,10){\makebox(0,0){\elt{Q}}}  %  bottom left
\put(10,35){\makebox(0,0){$B$}}  %  top left
\put(47,35){\makebox(0,0){\free{X}}}  % top right
\put(47,10){\makebox(0,0){\free{S}}}  % bottom right

\put(13,35){\vector(1,0){25}}  %     top
\put(16,10){\vector(1,0){23}}  %     bottom
\put(10,32){\vector(0,-1){17}} %     left
\put(45,32){\vector(0,-1){17}} %     right

\put(8,23){\makebox(0,0)[r]{}} %           left
\put(47,23){\makebox(0,0)[l]{}}%   right
\put(27,37){\makebox(0,0)[b]{}}%           top
\put(27,8){\makebox(0,0)[t]{}} %           bottom

\end{picture}
\end{center} and use the universal property of the pullback $A$ to
induce a morphism $B \lra A$, and hence $X'' \lra X'$.

The morphism $B \lra \free{X}$ along the top is given by
    \[\elt{Q} \times \free{X'} \map{p_2} \free{X'} \map{\cl{F}p_2}
    \free{\cl{F}X} \map{\mu} \free{X}\]
where $p_1$ and $p_2$ denote the first and second projections
respectively.  The morphism $B \lra \elt{Q}$ on the left is given
by
    \[\elt{Q} \times \free{X'} \map{(1, \cl{F} p_1)} \elt{Q}
    \times \free{(\elt{Q})} \lra \elt{Q}\]
where the second morphism is composition in $Q$.  Then, by
definition of $X'$ and naturality of $\mu$, the above square
commutes, inducing a map
    \[B \lra A\]
and hence, on isomorphism classes, a map

\setlength{\unitlength}{1mm}

\begin{center}
\begin{picture}(40,25)(0,10)      %

\put(20,12){\makebox(0,0){$S$}}  %  bottom
\put(5,30){\makebox(0,0){$X''$}}  %  top left
\put(35,30){\makebox(0,0){$X'$}}  %  top right

\put(8,30){\vector(1,0){24}}  %  top
\put(5,27){\vector(1,-1){13}}  %  left
\put(34,27){\vector(-1,-1){13}} %  right

\put(30,20){\makebox(0,0)[l]{$$}} %  right
\put(9,20){\makebox(0,0)[r]{$$}}%  left
\put(20,31){\makebox(0,0)[b]{$\mu_{(X,f)}$}}%  top

\end{picture}
\end{center}

\noindent in $\set/S$ as required.

Informally, $(X,f)$ is a system for labelling $Q$-objects, and
$T(X,f) = (X',f')$ gives source-labelled $Q$-arrows.  A typical
element of $X'$ may be thought of as the isomorphism class of

\setlength{\unitlength}{0.4mm}
\begin{center}
\begin{picture}(70,80)(15,40)      %

\put(20,90){\line(0,1){20}}      %
\put(40,90){\line(0,1){20}}      %
\put(80,90){\line(0,1){20}}      %

\put(20,90){\line(1,0){60}}      %
\put(20,90){\line(1,-1){30}}     %
\put(80,90){\line(-1,-1){30}}    %
\put(50,40){\line(0,1){20}}      %

\put(20,114){\makebox[0pt]{$x_1$}}   %
\put(40,114){\makebox[0pt]{$x_2$}}   %
\put(60,114){\makebox[0pt]{$\cdots$}}   %
\put(80,114){\makebox[0pt]{$x_m$}}   %

\put(50,76){\makebox(0,0){$\alpha$}}      % label in the middle
\put(50,36){\makebox(0,0)[t]{$$}}   % label at the bottom

\end{picture}
\end{center}

\noindent where $\alpha \in \elt{Q}$ and $s(\alpha) \cong
(\xv{fx})$.  Then $f'$ takes this element to $[t(\alpha)]$. So a
typical element $\theta$ of $T^2(X,f) = (X'',f'')$ is the
isomorphism class of

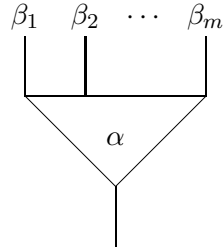
\begin{center}
\setlength{\unitlength}{0.4mm}
\begin{picture}(70,80)(15,40)    %

\put(20,90){\line(0,1){20}}      %
\put(40,90){\line(0,1){20}}      %
\put(80,90){\line(0,1){20}}      %

\put(20,90){\line(1,0){60}}      %
\put(20,90){\line(1,-1){30}}     %
\put(80,90){\line(-1,-1){30}}    %
\put(50,40){\line(0,1){20}}      %

\put(20,114){\makebox[0pt]{$\beta_1$}}   %
\put(40,114){\makebox[0pt]{$\beta_2$}}   %
\put(60,114){\makebox[0pt]{$\cdots$}}   %
\put(80,114){\makebox[0pt]{$\beta_m$}}   %

\put(50,76){\makebox(0,0){$\alpha$}}      % label in the middle
\put(50,36){\makebox(0,0)[t]{$$}}   % label at the bottom

\end{picture}
\end{center}where $\beta_i \in X'$ and $s(\alpha) \cong
(f'(\beta_1), \ldots, f'(\beta_m))$.  Writing $\beta_i$ as the
isomorphism class of

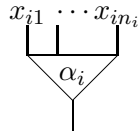
\begin{center}
\setlength{\unitlength}{0.2mm}
\begin{picture}(70,100)(15,40)      %

\put(20,90){\line(0,1){20}}      %
\put(40,90){\line(0,1){20}}      %
\put(80,90){\line(0,1){20}}      %

\put(20,90){\line(1,0){60}}      %
\put(20,90){\line(1,-1){30}}     %
\put(80,90){\line(-1,-1){30}}    %
\put(50,40){\line(0,1){20}}      %

\put(20,114){\makebox[0pt]{$x_{i1}$}}   %
\put(52,114){\makebox[0pt]{$\cdots$}}   %
\put(80,114){\makebox[0pt]{$x_{in_i}$}}   %

\put(50,76){\makebox(0,0){$\alpha_i$}}      % label in the middle
\put(50,36){\makebox(0,0)[t]{$$}}   % label at the bottom

\end{picture}
\end{center}we can draw $\theta$ as (the isomorphism class of)

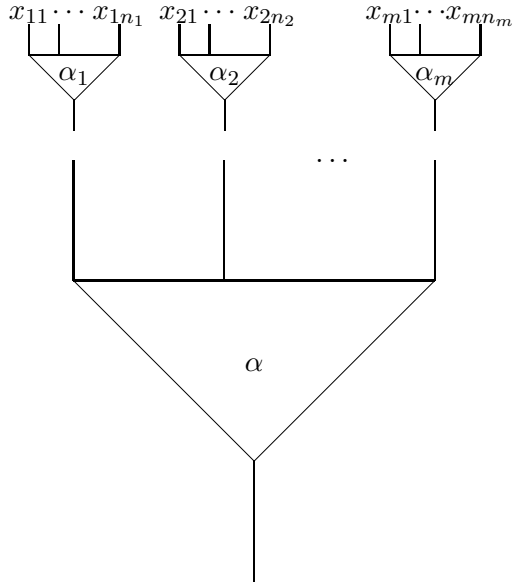
\begin{center}
\setlength{\unitlength}{0.8mm} %
\begin{picture}(70,100)(28,50)

\put(10,10){
\begin{picture}(90,140)      %
\put(20,90){\line(0,1){20}}      %
\put(45,90){\line(0,1){20}}      %
\put(63,110){\makebox(0,0){$\ldots$}}  %
\put(80,90){\line(0,1){20}}      %
\put(20,90){\line(1,0){60}}      %
\put(20,90){\line(1,-1){30}}     %
\put(80,90){\line(-1,-1){30}}    %
\put(50,40){\line(0,1){20}}      %
\put(50,76){\makebox(0,0){$\alpha$}}      % label in the middle
\end{picture}}  %

\put(16,115){ \setlength{\unitlength}{0.2mm}
\begin{picture}(90,140)      %
\put(20,90){\line(0,1){20}}      %
\put(40,90){\line(0,1){20}}      %
\put(80,90){\line(0,1){20}}      %
\put(20,90){\line(1,0){60}}      %
\put(20,90){\line(1,-1){30}}     %
\put(80,90){\line(-1,-1){30}}    %
\put(50,40){\line(0,1){20}}      %
\put(20,114){\makebox[0pt]{$x_{11}$}}   %
\put(48,114){\makebox[0pt]{$\cdots$}}   %
\put(80,114){\makebox[0pt]{$x_{1n_1}$}}   %
\put(50,76){\makebox(0,0){$\alpha_1$}}      % label in the middle
\put(50,36){\makebox(0,0)[t]{$$}}   % label at the bottom
\end{picture}}

\put(41,115){ \setlength{\unitlength}{0.2mm}
\begin{picture}(90,140)      %
\put(20,90){\line(0,1){20}}      %
\put(40,90){\line(0,1){20}}      %
\put(80,90){\line(0,1){20}}      %
\put(20,90){\line(1,0){60}}      %
\put(20,90){\line(1,-1){30}}     %
\put(80,90){\line(-1,-1){30}}    %
\put(50,40){\line(0,1){20}}      %
\put(20,114){\makebox[0pt]{$x_{21}$}}   %
\put(48,114){\makebox[0pt]{$\cdots$}}   %
\put(80,114){\makebox[0pt]{$x_{2n_2}$}}   %
\put(50,76){\makebox(0,0){$\alpha_2$}}      % label in the middle
\put(50,36){\makebox(0,0)[t]{$$}}   % label at the bottom
\end{picture}}

\put(76,115){ \setlength{\unitlength}{0.2mm}
\begin{picture}(90,140)      %
\put(20,90){\line(0,1){20}}      %
\put(40,90){\line(0,1){20}}      %
\put(80,90){\line(0,1){20}}      %
\put(20,90){\line(1,0){60}}      %
\put(20,90){\line(1,-1){30}}     %
\put(80,90){\line(-1,-1){30}}    %
\put(50,40){\line(0,1){20}}      %
\put(20,114){\makebox[0pt]{$x_{m1}$}}   %
\put(48,114){\makebox[0pt]{$\cdots$}}   %
\put(80,114){\makebox[0pt]{$x_{mn_m}$}}   %
\put(50,76){\makebox(0,0){$\alpha_m$}}      % label in the middle
\put(50,36){\makebox(0,0)[t]{$$}}   % label at the bottom
\end{picture}}

\end{picture}
\end{center}where $\alpha, \alpha_1, \ldots, \alpha_m \in \elt{Q}$
and $s(\alpha) \cong (t(\alpha_1), \ldots t(\alpha_m))$. So, via
the relevant object-isomorphisms, we may compose the underlying
$Q$-arrows to give $\alpha'$, say, which is defined up to
isomorphism.  We then concatenate the $X$-labels (via the
multiplication for \cl{F}) to give

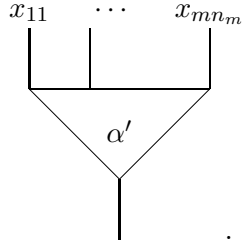
\begin{center}
\setlength{\unitlength}{0.4mm}
\begin{picture}(70,80)(15,40)       %

\put(20,90){\line(0,1){20}}      %
\put(40,90){\line(0,1){20}}      %
\put(80,90){\line(0,1){20}}      %

\put(20,90){\line(1,0){60}}      %
\put(20,90){\line(1,-1){30}}     %
\put(80,90){\line(-1,-1){30}}    %
\put(50,40){\line(0,1){20}}      %

\put(20,114){\makebox[0pt]{$x_{11}$}}   %
\put(48,114){\makebox[0pt]{$\cdots$}}   %
\put(80,114){\makebox[0pt]{$x_{mn_m}$}}   %

\put(50,76){\makebox(0,0){$\alpha'$}}      % label in the middle
\put(50,36){\makebox(0,0)[t]{$$}}   % label at the bottom

\end{picture}.
\end{center}Finally, we take the isomorphism class of this to
give $\mu_{(X,f)}(\theta) \in X'$, and $f''(\mu_{(X,f)}(\theta)) =
[t(\alpha')] = [t(\alpha)] \in S$.

It follows that $\mu$ defined in this way is a cartesian natural
transformation.

\begin{itemize}
\item $T$ preserves pullbacks
\end{itemize}

First observe that a commutative square in $\set/S$ is a pullback
if and only if applying the forgetful functor $\set/S
\longrightarrow \set$ gives a pullback in \set. Then $T$ preserves
pullbacks since \cl{F} preserves pullbacks.

\bigskip

So $T_Q=(T,\eta,\mu)$ is a cartesian monad on $\cl{E}_Q=\cl{E}$
and we may define $\zeta(Q)=(\cl{E}_Q,T_Q)$.

\bigskip

We now define the action of $\zeta$ on morphisms.  Let \[F:Q
\longrightarrow R\] be a morphism of tidy symmetric
multicategories. We construct a cartesian monad opfunctor
    \[(U_F,\phi_F) : (\cl{E}_Q,T_Q) \longrightarrow
    (\cl{E}_R,T_R)\]
that is \begin{itemize}
\item a functor $U=U_F: \set/S_Q \longrightarrow \set/S_R$
preserving pullbacks
\item a cartesian natural transformation $\phi=\phi_F: UT_Q
\longrightarrow T_RU$
\end{itemize} satisfying certain axioms.

We define $U$ as follows. On objects, we have a functor
    \[F:o(Q) \longrightarrow o(R)\]
giving a morphism on isomorphism classes
    \[\bar{F}:S_Q \longrightarrow S_R.\]
This induces a functor
    \[\set/S_Q \longrightarrow \set/S_R\]
by composition with $\bar{F}$, which clearly preserves pullbacks;
we define $U$ to be this functor.

We now construct the components of $\phi$.  Given $(X,f) \in
\set/S_Q$ write
\[\begin{array}{c}
    T_Q(X,f) = (X^Q,f^Q)\\
    \makebox(0,0)[br]{and \ \ }X^Q \simeq \xx{Q}{S_Q}{X}.
\end{array}\]
We seek
    \[\phi_{(X,f)} : (X^Q,\bar{F}\circ f^Q) \longrightarrow (X^R,
    (\bar{F} \circ f)^R) \in \set/S_R\]
that is, a morphism $X^Q \longrightarrow X^R$ such that the
outside of the following diagram commutes

\begin{center}
\setlength{\unitlength}{1mm}
\begin{picture}(60,115)      %

\put(10,10){\makebox(0,0){$S_Q$}}  %  bottom left
\put(10,35){\makebox(0,0){$o(Q)$}}  %  top left
\put(45,35){\makebox(0,0){$o(R)$}}  %  top right
\put(45,10){\makebox(0,0){$S_R$}}  %  bottom right

\put(10,60){\makebox(0,0){$\elt{Q}$}}  %
\put(45,60){\makebox(0,0){$\elt{R}$}}  %
\put(6,85){\makebox(0,0){$\xx{Q}{S_Q}{X}$}}  %
\put(49,85){\makebox(0,0){$\xx{R}{S_R}{X}$}}  %
\put(10,110){\makebox(0,0){$X^Q$}}  %
\put(45,110){\makebox(0,0){$X^R$}}  %

\put(13,10){\vector(1,0){28}}  %  bottom
\put(15,35){\vector(1,0){24}}  %
\put(15,60){\vector(1,0){24}}  %
\put(21,85){\vector(1,0){10}}  %
%\put(13,110){\dashbox{0.3}(28,0){}} % top
\qbezier[40](13,110)(28,110)(41,110)%
\put(40,110){\vector(1,0){2}}   %

\put(10,32){\vector(0,-1){19}} %  left
\put(45,32){\vector(0,-1){19}} %  right
\put(10,57){\vector(0,-1){19}} %  left
\put(45,57){\vector(0,-1){19}} %  right
\put(10,82){\vector(0,-1){19}} %  left
\put(45,82){\vector(0,-1){19}} %  right
\put(10,107){\vector(0,-1){19}} %  left
\put(45,107){\vector(0,-1){19}} %  right

\put(8,23){\makebox(0,0)[r]{$\sim$}} %  left
\put(47,23){\makebox(0,0)[l]{$\sim$}}%  right
\put(8,48){\makebox(0,0)[r]{$t$}} %  left
\put(47,48){\makebox(0,0)[l]{$t$}}%  right
\put(8,98){\makebox(0,0)[r]{$\sim$}} %  left
\put(47,98){\makebox(0,0)[l]{$\sim$}}%  right

\put(27,8){\makebox(0,0)[t]{$\bar{F}$}} %  bottom
\put(27,37){\makebox(0,0)[b]{$F$}}%  top
\put(27,62){\makebox(0,0)[b]{$F$}}%
\put(27,87){\makebox(0,0)[b]{$(F,1)$}}%

\end{picture}.
\end{center}The map $X^Q \longrightarrow X^R$ is induced by $(F,1)$ on
isomorphism classes as shown in the diagram, since the pullback
    \[\xx{R}{S_R}{X}\]
is along the morphism $\bar{F} \circ f$. We define $\phi_{(X,f)}$
to be this map. Observe that all squares in the diagram commute,
so $\phi_{(X,f)}$ is a morphism in $\set/S_R$ as required.

%11

We now check that these components satisfy naturality.  Given any
morphism $h:(X,f) \longrightarrow (Y,g) \in \set/S_Q$, we have the
following diagram

\begin{center}
\setlength{\unitlength}{1mm}
\begin{picture}(80,45)(-4,0)       %

\put(10,10){\makebox(0,0){$\xx{Q}{S_Q}{Y}$}}  %  bottom left
\put(10,35){\makebox(0,0){$\xx{Q}{S_Q}{X}$}}  %  top left
\put(60,35){\makebox(0,0){$\xx{R}{S_R}{X}$}}  %  top right
\put(60,10){\makebox(0,0){$\xx{R}{S_R}{Y}$}}  %  bottom right

\put(25,35){\vector(1,0){18}}  %  top
\put(25,10){\vector(1,0){18}}  %  bottom
\put(9,32){\vector(0,-1){19}} %  left
\put(60,32){\vector(0,-1){19}} %  right

\put(8,23){\makebox(0,0)[r]{$(1,\cl{F}h)$}} %  left
\put(62,23){\makebox(0,0)[l]{$(1,\cl{F}h)$}}%  right
\put(33,37){\makebox(0,0)[b]{$(F,1)$}}%  top
\put(33,8){\makebox(0,0)[t]{$(F,1)$}} %  bottom

\end{picture}
\end{center}Considering this componentwise, it clearly commutes and is a
pullback.  The result on isomorphism classes follows.

Finally, by functoriality of $F$, $(U,\phi)$ satisfies the axioms
for a monad opfunctor.  So $(U,\phi)$ is a cartesian monad
opfunctor and the construction is clearly functorial.  This
completes the definition of $\zeta$.

\bigskip

We observe immediately that the construction of $(\cl{E}_Q, T_Q)$
uses only the isomorphism classes of objects and arrows of $Q$. So
    \[(\cl{E}_{Q_1}, T_{Q_1}) \cong (\cl{E}_{Q_2}, T_{Q_2}) \iff
    Q_1 \simeq Q_2.\]

Recall (\ref{tacsym}) that we expect that a symmetric
multicategory $Q$ may be given as a monad in a certain bicategory,
in which case the identities are given by the unit, and
composition laws by multiplication.  In this abstract framework
there should be a morphism from the underlying bicategory to the
$2$-category \cat{Cat}, taking the monad $Q$ to the monad $T_Q$,
but this is somewhat beyond the scope of this work.

\section{The theory of opetopes} \label{opes}

In this section we give the analogous constructions of opetopes in
each theory, and show in what sense they are equivalent.  That is,
we show that the respective categories of $k$-opetopes are
equivalent.

We first discuss the process by which $(k+1)$-cells are
constructed from $k$-cells.  Recall that, in \cite{bd1}, the `slice'
construction is used, giving for any \sm\ $Q$ the slice \mcat\
$Q^+$.  In \cite{lei1} the `free $(\cl{E},T)$-operad' construction is
used, giving, for any `suitable' monad \st, the free
$(\cl{E},T)$-operad monad $(\cl{E}',T')$.

\subsection{Slicing a symmetric multicategory}
\label{tacsymslice}

Let $Q$ be a symmetric multicategory with a category ${\mathbb C}$ of
objects, so $Q$ may be considered as a functor $Q:\cl{F}{\mathbb
C}^{\mbox{\scriptsize op}} \times {\mathbb C} \longrightarrow
\mbox{\cat{Set}}$ with certain extra structure. Recall (\cite{che7})
that the slice multicategory $Q^+$ is given by:

\begin{itemize}
\item Objects: put $o(Q^+) = \mbox{elt}(Q)$
\end{itemize}

So the category $o(Q^+)$ has as objects pairs $(p,g)$ with $p \in
\cl{F}{\mathbb C}^{\mbox{\scriptsize op}} \times {\mathbb C}$ and
$g \in Q(p)$; a morphism $\alpha:(p,g) \longrightarrow (p',g')$ is
an arrow
        $\alpha:p \longrightarrow p' \in \cl{F}{\mathbb C}
        ^{\mbox{\scriptsize op}} \times {\mathbb C}$
such that
    \[\begin{array}{rccc}
    Q(\alpha) : & Q(p) & \lra & Q(p')
    \\
    & g & \longmapsto & g' \  \\
    \end{array}\]

Then, given any arrow \[g \in Q(x_1, \ldots x_m; x)\] we have an
arrow $\alpha(g) = g' \in Q(y_1, \ldots, y_m; y)$ given by
    \[g' = (\iota(f) \circ g \circ (\iota(f_1), \ldots,
    \iota(f_m))\sigma)\]
(see \cite{che7}).

\begin{itemize}
\item Arrows: $Q^+(f_1, \ldots, f_n;f)$ is given by the set of
`configurations' for composing $f_1, \ldots, f_n$ as arrows of
$Q$, to yield $f$.
\end{itemize}

Writing $f_i \in Q(x_{i1}, \ldots x_{im_i}; x_i)$ for $1 \leq i
\leq n$, such a configuration is given by $(T,\rho, \tau)$ where

\begin{enumerate}

\item $T$ is a planar tree with $n$ nodes.  Each node is labelled
by one of the $f_i$, and each edge is labelled by an
object-morphism of $Q$ in such a way that the (unique) node
labelled by $f_i$ has precisely $m_i$ edges going in from above,
labelled by $a_{i1}, \ldots, a_{im_i} \in \mbox{arr}({\mathbb
C})$, and the edge coming out is labelled $a_i \in a({\mathbb
C})$, where $\mbox{cod}(a_{ij}) = x_{ij}$ and $\mbox{dom}(a_i) =
x_i$.

\item $\rho \in {\mathbf S}_k$ where $k$ is the number of leaves
of $T$.

\item $\tau:\{\mbox{nodes of } T\} \longrightarrow [n]=\{1, \ldots, n\}$ is a
bijection such that the node $N$ is labelled by $f_{\tau(N)}$.
(This specification is necessary to allow for the possibility $f_i
= f_j,\ i \neq j$.)

\end{enumerate} Note that $(T,\rho)$ may be considered as a `combed tree',
that is, a planar tree with a `twisting' of branches at the top
given by $\rho$.

The arrow resulting from this composition is given by composing
the $f_i$ according to their positions in $T$, with the $a_{ij}$
acting as arrows $\iota(a_{ij})$ of $Q$, and then applying $\rho$
according to the symmetric action on $Q$. This construction
uniquely determines an arrow $(T,\rho,\tau) \in Q^+(f_1, \ldots,
f_n;f)$.

\begin{itemize}
\item Composition
\end{itemize}

When it can be defined, $(T_1,\rho_1,\tau_1)
    \circ_m (T_2,\rho_2,\tau_2) = (T,\rho,\tau)$ is given by

\begin{enumerate}

\item $(T,\rho)$ is the combed tree obtained by replacing the node
${\tau_1}^{-1}(m)$ by the tree $(T_2,\rho_2)$, composing the edge
labels as morphisms of ${\mathbb C}$, and then `combing' the tree
so that all twists are at the top.

\item $\tau$ is the bijection which inserts the source of $T_2$
into that of $T_1$ at the $m$th place.

\end{enumerate}

\begin{itemize}

\item Identities: given an object-morphism
    \[\alpha=(\sigma, f_1, \ldots, f_m;f) : g \longrightarrow g',\]
$\iota(\alpha) \in Q^+(g;g')$ is given by a tree with one node,
labelled by $g$, twist $\sigma$, and edges labelled by the $f_i$
and $f$ as in the example above.

\item Symmetric action: $(T,\rho,\tau)\sigma = (T,\rho,\sigma^{-1}\tau)$

\end{itemize}

\noindent This is easily seen to satisfy the axioms for a
symmetric multicategory.

Recall (\cite{che7}) that $Q^+$ is always freely symmetric, and that
it is tidy if $Q$ is tidy.

\subsection{Slicing a \tmcat}
\label{tacslicing}

In \cite{lei1} the `free $(\cl{E},T)$-operad' construction is used
to construct $(k+1)$-cells from $k$-cells; this gives, for any
{\em suitable} monad \st, the `free $(\cl{E},T)$-operad' monad
$(\cl{E},T)' = (\cl{E}',T')$.  In order to compare this
construction with the Baez-Dolan slice, we examine the monad
$\zeta(Q)'$. First we must show that ${\zeta(Q)}'$ can actually be
constructed, that is, that ${\zeta(Q)}=\sqtq$ is a suitable monad.

First recall (\cite{lei1}) that a cartesian monad $(\cl{E}, T)$ is
{\em suitable} if it satisfies:
\begin{enumerate}
\item $\cl{E}$ has disjoint finite coproducts which are stable
under pullback
\item $\cl{E}$ has colimits of nested sequences; these commute
with pullbacks and have monic coprojections
\item $T$ preserves colimits of nested sequences.
\end{enumerate}
Here a nested sequence is a string of composable monics.

\begin{proposition} \label{opespropw} Let $Q$ be a tidy \sm.  Then
\sqtq\ is a suitable monad.
\end{proposition}

\begin{prf} Certainly $\cl{E}_Q$ is a suitable category, and we
have already shown that \sqtq\ is cartesian. So it remains to show
that $T_Q$ preserves colimits of nested sequences.

First observe that a morphism $h$ in $\set/S$ is monic if and only
if $h$ is monic as a morphism in \set, that is, injective. Given a
nested sequence
    \[(A_0,f_0) \stackrel{i_0}{\monic} (A_1,f_1)
    \stackrel{i_1}{\monic} (A_2,f_2) \cdots \in \set/S\]
we have a nested sequence
    \[A_0 \stackrel{i_0}{\monic} A_1
    \stackrel{i_1}{\monic} A_2 \cdots \in \set.\]
Since \set\ is suitable, this nested sequence has a colimit $A$
whose coprojections are monics.  Then the morphisms $f_0,f_1,
\ldots$ define a cone with vertex $S$, inducing a unique morphism
$A \map{f} S$ making everything commute; $(A,f)$ is then a colimit
for the nested sequence in $\set/S$. So $(A,f)$ is a colimit for
the nested sequence in $\set/S$ exactly when $A$ is a colimit for
the nested sequence in \set.

Having made these observations, it is easy to check that $T_Q$
preserves such colimits. \end{prf}

\subsection{Comparison of slice}
\label{slicesymcart}

In this section we compare the slice constructions and make
precise the sense in which they correspond to one another. Recall
(Section~\ref{taczeta}) that we have defined a functor
    \[\cat{TidySymMulticat} \map{\zeta} \cat{CartMonad}.\]
We now show that this functor `commutes' with slicing, up to
isomorphism (for $\zeta$).

Since $\zeta(Q)=\sqtq$ is suitable (Proposition~\ref{opespropw}), we can
form $\zeta(Q)'=\tqprime$, the free \tqopd\ monad. Also, $Q^+$ is
tidy since $Q$ is tidy (see \cite{che7}), so we can
form the monad $\zeta(Q^+) = \tqplus$. For the comparison, we have
the following result.

\begin{proposition} \label{tacpropx} Let Q be a tidy \sm.  Then
    \[\zeta(Q)' \cong \zeta(Q^+)\]
that is
    \[\tqprime \cong \tqplus\]
in the category \cat{CartMonad}.
\end{proposition}

\sunit This proof is somewhat technical and we defer it to
Appendix~\ref{apa}.  Informally, the idea is as follows. $T_{Q^+}$
takes a set $A$ of `labels for arrows of Q' and returns the set
$A_2$ of configurations for composing labelled arrows according to
their underlying arrows.  On the other hand, $T_Q'$ takes a
diagram of the form
\begin{center}
\begin{picture}(55,45)(8,2)  %
\tomspan{S}{T_Q}{S}{1}{A}{}{}{}
\end{picture}
\end{center}
and forms the free $(\cl{E}_Q, T_Q)$ multicategory on it, with
underlying graph
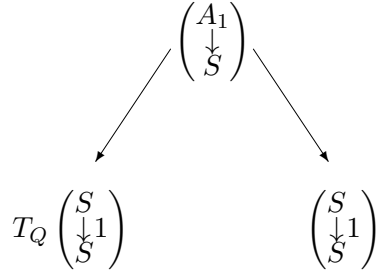
\begin{center}
\begin{picture}(55,45)(8,2)  %
\tomspan{S}{T_Q}{S}{1}{A_1}{}{}{}
\end{picture}.
\end{center}

So $T_Q'$ gives the set $A_1$ of all formal composites of arrows
labelled in $A$ according to the structure of $T_Q$, which is
precisely the set of configurations as above.

\bigskip
Recall that
    \[\zeta(Q_1) \cong \zeta(Q_2) \iff Q_1 \simeq Q_2.\]
We immediately deduce the following result, comparing all three
processes of slicing.

\begin{corollary} Let $M$ be a generalised multicategory.  Then
    \[\zeta \xi (M_+) \cong \zeta (\xi(M)^+) \cong \zeta \xi
    (M)'.\]
\end{corollary}

\label{tacopetopes}

We are now ready to compare the different constructions of opetopes,
applying the results we have already established. In each case,
opetopes are constructed by iterating the slicing process. Note that
the `opetopes' defined in \cite{lei1} are not {\em a priori} the same
as those defined in \cite{bd1}; when a distinction is required we
refer to the former as `Leinster opetopes'.

\subsection{Opetopes}
\label{tacopes}

We recall the definition of opetopes from \cite{bd1} and \cite{che7}. 
For any symmetric multicategory $Q$ we write
    \[Q^{k+} = \left\{\begin{array} {l@{\extracolsep{2em}}l}
    Q & k=0 \\
    {(Q^{(k-1)+})}^+ & k \ge 1\end{array} \right. \]
Let $I$ be the symmetric multicategory with precisely one object,
precisely one (identity) object-morphism, and precisely one
(identity) arrow.  A {\em $k$-dimensional opetope}, or simply {\em
$k$-opetope}, is defined in \cite{bd1} to be an object of
$I^{k+}$. We write $\bb{C}_k = o(I^{k+})$, the category of
$k$-opetopes.

\subsection{Leinster opetopes}
\label{leiopes}

In \cite{lei1}, $k$-opetopes are defined by a sequence
$(\cat{Set}/S_k, T_k)$ of cartesian monads given by iterating the
slice as follows.

For any cartesian monad (\cl{E},T) write
\[(\cl{E},T)^{k'} = \left\{\begin{array} {l@{\extracolsep{2em}}l}
    (\cl{E},T) & k=0 \\
    {((\cl{E},T)^{(k-1)'})}' & k \ge 1 \end{array} \right.\]
Put $(\cl{E}_0,T_0)=(\cat{Set},\id)$ and for $k \ge 1$ put
$(\cl{E}_k,T_k)=(\cat{Set},\id)^{k'}$.  It follows that for each
$k$, $(\cl{E}_k,T_k)$ is of the form $(\cat{Set}/S_k, T_k)$ where
$S_0=1$ and $S_{k+1}$ is given by
 \[\vect{S_{k+1}}{}{S_k} = T_k \vect{S_k}{1}{S_k}\]
Then {\em Leinster $k$-opetopes} are defined to be the elements of
$S_k$; as above, we will regard $S_k$ as a discrete category.

\subsection{Comparisons of opetopes} \label{opecomp}

We now compare opetopes and Leinster opetopes.

\begin{proposition} \label{tacpropy} For each $k \geq 0$
    \[\zeta(I^{k+}) \cong ({\mbox{\upshape{\bfseries{Set}}}},
    \id)^{k'} = (\cat{Set}/S_k, T_k).\]
\end{proposition}

\begin{prf} By induction. For $k=0$ we need to show
    \[(\cl{E}_{I^{k+}}, T_{I^{k+}}) \cong (\set, \id).\]
Now $\cl{E}_I = \set/S_I$ where $S_I \simeq o(I) = 1$.  So
$\cl{E}_I \cong \set/1 \cong \set$. Given any $\vect{X}{!}{1} \in
\set/1$,  $T_I \vect{X}{}{1}$ is equivalent to the pullback

\setlength{\unitlength}{1mm}
\begin{center}
\begin{picture}(45,35)(5,5)    %

\put(10,10){\makebox(0,0){$\elt{I}$}}  %  bottom left
\put(10,35){\makebox(0,0){$\cdot$}}  %  top left
\put(45,35){\makebox(0,0){$\free{X}$}}  %  top right
\put(45,10){\makebox(0,0){$\free{1}$}}  %  bottom right

\put(13,35){\vector(1,0){26}}  %  top
\put(14,10){\vector(1,0){25}}  %  bottom
\put(10,32){\vector(0,-1){19}} %  left
\put(45,32){\vector(0,-1){19}} %  right

\put(8,23){\makebox(0,0)[r]{$$}} %  left
\put(47,23){\makebox(0,0)[l]{$$}}%  right
\put(27,37){\makebox(0,0)[b]{$$}}%  top
\put(27,8){\makebox(0,0)[t]{$t$}} %  bottom

\end{picture}.
\end{center} But $I$ has only one arrow, which is unary (the identity), so
    \[T_I \vect{X}{}{1} \cong \vect{X}{}{1}\]
and
    \[(\cl{E}_I, T_I) \cong (\set,\id)\]
as required.

Now suppose $\zeta(I^{(k-1)+}) \cong (\set, \id)^{(k-1)'}$.  Then
by Proposition~\ref{tacpropx} we have
    \[\zeta(I^{k+}) \cong \zeta(I^{(k-1)+})' \cong (\set,
    \id)^{k'}\]
so by induction the result is true for all $k \geq 0$.
\end{prf}

Then on objects, the above equivalence gives the following result.

\begin{corollary} For each $k \geq 0$
    \[{\mathbb C}_k \simeq S_k.\]
\end{corollary}

Recall (\cite{che7}) that we also have for each $k$ a (discrete)
category $P_k$ of `multitopes', the analogous notion as defined in
\cite{hmp1}; in \cite{che7} we prove that, for each $k \geq 0$, $P_k
\simeq {\mathbb C}_k$.  So we immediately have the following result,
comparing all three theories:

\begin{corollary} For each $k \geq 0$
    \[P_k \simeq {\mathbb C}_k \simeq S_k.\]
\end{corollary}

This result shows that multitopes, opetopes, and Leinster opetopes are
the same, up to isomorphism.

We eventually aim to define a category \cat{Opetope} of opetopes
of all dimensions, whose morphisms are `face maps' of opetopes;
this is the subject of \cite{che9}.

%%%%%%%%%%%%%%%%%%%

\appendix

\section{Proof of Proposition~\ref{tacpropx}}
\label{apa}

We now give the proof of Proposition~\ref{tacpropx} deferred from 
Section~\ref{slicesymcart}. 

\bigskip
\noindent {\bfseries Proposition \ref{tacpropx}.  } Let Q be a 
tidy \sm.  Then 
    \[\zeta(Q)' \cong \zeta(Q^+)\]
that is 
    \[\tqprime \cong \tqplus\]
in the category \cat{CartMonad}.

\begin{prf}  First we show that ${\cl{E}_Q}' \cong \cl{E}_{Q^+}$. Now 
$\cl{E}_{Q^+} = \set/S_{Q^+}$ where $S_{Q^+} \cong 
o(Q^+)=\elt{Q}$, and ${\cl{E}_Q}' = \set/{S_Q}'$ where 
\[\vect{{S_Q}'}{}{S_Q} = T_Q \vect{S_Q}{1}{S_Q}.\] So ${S_Q}'$ is
equivalent to the pullback 

\begin{center}
\begin{picture}(45,35)(5,5)     %

\put(10,10){\makebox(0,0){$\elt{Q}$}}  %  bottom left
\put(10,35){\makebox(0,0){$\cdot$}}  %  top left
\put(47,35){\makebox(0,0){$\free{S_Q}$}}  %  top right
\put(47,10){\makebox(0,0){$\free{S_Q}$}}  %  bottom right

\put(13,35){\vector(1,0){26}}  %  top
\put(14,10){\vector(1,0){25}}  %  bottom
\put(10,32){\vector(0,-1){19}} %  left
\put(45,32){\vector(0,-1){19}} %  right

\put(8,23){\makebox(0,0)[r]{$$}} %  left
\put(47,23){\makebox(0,0)[l]{$1$}}%  right
\put(27,37){\makebox(0,0)[b]{$$}}%  top
\put(27,8){\makebox(0,0)[t]{$$}} %  bottom

\end{picture}
\end{center} so ${S_Q}' \simeq \elt{Q}$, giving ${S_Q}' \cong S_{Q^+}$.
So we have ${\cl{E}_Q}' \cong \cl{E}_{Q^+}$.  By abuse of 
notation, we write elements of both these categories as sets over 
$S'$, since confusion is unlikely. 

Consider $(A,f) = (A \map{f} S') \in {\cl{E}_Q}' \cong 
\cl{E}_{Q^+}$.  Write $\tpr(A,f) = (A_1,f_1)$ and $\tpl(A,f) = 
(A_2,f_2)$. We show $(A_1,f_1) \cong (A_2,f_2)$. To construct 
$A_2$, first form the pullback 

\setlength{\unitlength}{1mm} 
\begin{center}
\begin{picture}(70,35)(5,5)  %

\put(10,10){\makebox(0,0){$\elt{Q^+}$}}  %  bottom left
\put(10,35){\makebox(0,0){$\cdot$}}  %  top left
\put(72,35){\makebox(0,0){$\free{A}$}}  %  top right
\put(40,10){\makebox(0,0){$\free{(\elt{Q})}$}}  %  bottom middle
\put(72,10){\makebox(0,0){$\free{S'}$}}  %  bottom right

\put(13,35){\vector(1,0){52}}  %  top
\put(15,10){\vector(1,0){14}}  %  bottom left
\put(49,10){\vector(1,0){16}}  %  bottom right
\put(10,32){\vector(0,-1){19}} %  left
\put(70,32){\vector(0,-1){19}} %  right

\put(8,23){\makebox(0,0)[r]{$$}} %  left
\put(72,23){\makebox(0,0)[l]{\free{f}}}%  right
\put(27,37){\makebox(0,0)[b]{$$}}%  top
\put(22,8){\makebox(0,0)[t]{$s$}} %  bottom left
\put(56,8){\makebox(0,0)[t]{$\sim$}} %  bottom right

\end{picture}.
\end{center} Then $A_2 \simeq \xx{Q^+}{S'}{A}$, and $f_2$ is given by the
composite 
    \[ A_2 \simeq \xx{Q^+}{S'}{A} \longrightarrow \elt{Q^+} \map{t_{Q^+}}
    \elt{Q} \map{\sim} S'\]
where $t_{Q^+}$ is the target map of $Q^+$. 

Informally, since we are here considering $S' \simeq o(Q^+) = 
\elt(Q)$, the object $(A \map{f} S')$ may be thought of as a set 
of labels for arrows of $Q$. Then $A_2$ is the set of %isomorphism classes of %
all possible source-labelled arrows of $Q^+$. Since an arrow of 
$Q^+$ is given by a tree with nodes corresponding to arrows of 
$Q$, an element of $A_2$ may be thought of as %an isomorphism class of%
such a tree, with nodes labelled by compatible elements of $A$. 
Alternatively, it may be thought of as a configuration for 
composing labelled arrows of $Q$ via object-isomorphisms, where 
composition is according to the underlying arrows only. $f_2$ acts 
by composing the underlying arrows of $Q$ and then taking 
isomorphism classes. 

We now turn our attention to the action of \tpr.  (For full 
details of the free multicategory construction we refer the reader 
to \cite{lei4}.)  For convenience we write $T_Q = T$ and $S_Q=S$, 
so we need to form 
    \[(T, \set/S)' = (T', S').\]
To construct $A_1$, we form the free multicategory on the 
following graph: 

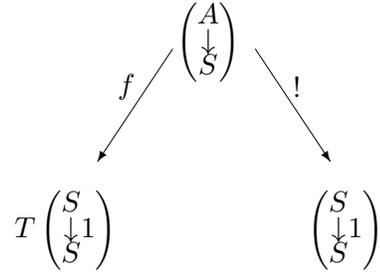
\length{1mm} 
\begin{center}
\begin{picture}(55,45)(8,2)  %
\tomspan{S}{T}{S}{1}{A}{}{f}{!} 
\end{picture}.
\end{center} 
Recall we have
    \[T\vect{S}{}{S} = \vect{S'}{}{S}\]
and the map $A \lra S$ is the composite $A \map{f} S' \lra S$. The 
graph underlying the free operad is then 

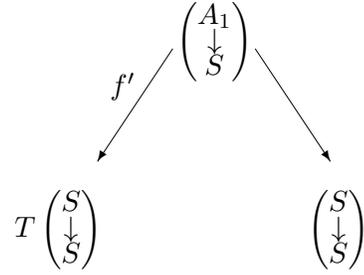
\begin{center}
\begin{picture}(55,45)(8,2)  %
\tomspan{S}{T}{S}{}{A_1}{}{f'}{} 
\end{picture}.
\end{center}

\noindent The construction gives a sequence of graphs 

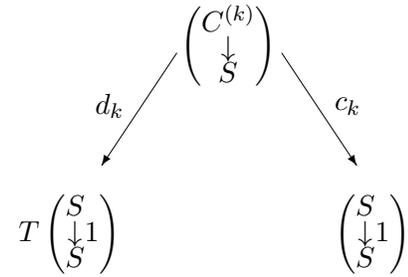
\setlength{\unitlength}{1mm} 
\begin{center}
\begin{picture}(55,45)(8,2)      %
\put(8,10){${T}\vect{S}{1}{S}$}  %  bottom left
\put(29.5,35){$\vect{C^{(k)}}{}{S}$}  %  top
\put(50,10){$\vect{S}{1}{S}$}  %  bottom right

\put(29,35){\vector(-2,-3){10}}  %
\put(43,35){\vector(2,-3){10}}  %
\put(22,28){\makebox(0,0)[r]{${d_k}$}} %  left
\put(50,28){\makebox(0,0)[l]{${c_k}$}}%  right

\end{picture}
\end{center} where $C^{(0)} = S$, $d_0 = \eta_T$ and
\[\vect{C^{(k+1)}}{}{S} = \vect{S}{1}{S} + \vect{A}{}{S} \circ
\vect{C^{(k)}}{}{S}.\]

%8

\noindent Here $\circ$ is composition in the bicategory of spans, 
so the composite \[\vect{A}{}{S} \circ \vect{C^{(k)}}{}{S}\] is 
given by the pullback 

\setlength{\unitlength}{1mm} 
\begin{center}
\begin{picture}(60,45) %

\put(10,10){\makebox(0,0){$T\vect{C^{(k)}}{}{S}$}}  %  bottom left
\put(10,35){\makebox(0,0){$\cdot$}}  %  top left
\put(47,35){\makebox(0,0){$\vect{A}{}{S}$}}  %  top right
\put(53,10){\makebox(0,0){$T\vect{S}{}{S}=\vect{S'}{}{S}$}}  %

\put(13,35){\vector(1,0){27}}  %  top
\put(18,10){\vector(1,0){20}}  %  bottom
\put(10,32){\vector(0,-1){14}} %  left
\put(47,28){\vector(0,-1){10}} %  right

\put(8,23){\makebox(0,0)[r]{$$}} %  left
\put(49,23){\makebox(0,0)[l]{$f$}}%  right
\put(27,37){\makebox(0,0)[b]{$$}}%  top
\put(27,8){\makebox(0,0)[t]{$Tc_k$}} %  bottom

\end{picture}
\end{center} and $d_{k+1}$ is given by the composite
    \[\vect{A}{}{S} \circ \vect{C^{(k)}}{}{S} \lra
        T\vect{C^{(k)}}{}{S} \map{Td_k}
        TT\vect{S}{}{S} \map{\mu_T} T\vect{S}{}{S} .\]
This construction gives a nested sequence $(C^{(k)},f^{(k)}) \in 
\set/S$ with $(C^{(0)},f^{(0)})=(S,1)$ and 
    \[C^{(k+1)} = S \ \amalg \ T(C^{(k)}) \times_{S'} A\]
where (by further abuse of notation) we write 
    \[T\vect{C^{(k)}}{}S{}=\vect{T(C^{(k)})}{}{S}.\] 
$f^{(k+1)}$ is given by $1 \ \amalg \ (T(C^{(k)}) \times_{S'} A 
\map{d_{k+1}} S' \lra S)$ and $\vect{A_1}{}{S}$ is then the 
colimit of this nested sequence. 

Informally, the sets $C^{(k)}$ may be thought of as $k$-fold 
formal composites (or composites of `depth' at most $k$).  The 
formula for $C^{(k)}$ says that a composite is either null or is a 
generating arrow composed with other composites.  We aim to show 
that these formal composites correspond to the formal composites 
given by the source-labelled arrows of $Q^+$. 

We show that $A_1 \cong A_2 \simeq \xx{Q^+}{S'}{A}$ as follows. 
For each $k$ we exhibit an embedding 
    \[g_k : C^{(k)} {\hookrightarrow} A_2\]
which makes the following diagram commute 

\setlength{\unitlength}{1mm} 
\begin{center}
\begin{picture}(38,25)(0,10)     %

\put(20,12){\makebox(0,0){$S'$}}  %  bottom
\put(5,30){\makebox(0,0){$C^{(k)}$}}  %  top left
\put(35,30){\makebox(0,0){$A_2$}}  %  top right

\put(9,30){\vector(1,0){22}}  %  top
\put(5,27){\vector(1,-1){13}}  %  left
\put(34,27){\vector(-1,-1){13}} %  right

\put(30,20){\makebox(0,0)[l]{$f_2$}} %  right
\put(9,20){\makebox(0,0)[r]{$d_k$}}%  left
\put(20,31){\makebox(0,0)[b]{$g_k$}}%  top

\end{picture}.
\end{center} Then the colimit induces the map required.

We proceed by induction.  Define $g_0:S \lra \xx{Q^+}{S'}{A}$ as 
follows. Let $[x] \in S$ denote the isomorphism class of $x \in 
o(Q)$. Given any $[x] \in S$, we have a nullary arrow $\alpha_x 
\in Q^+(\cdot\ ;1_x)$.  Recall that an arrow of $Q^+$ may be 
regarded as a tree with nodes corresponding to the source elements 
(which are themselves arrows of $Q$) and edges labelled by 
object-morphisms of $Q$.  Then $\alpha_x \in Q^+(\cdot\ ;1_x)$ is 
given by a tree with no nodes, that is, a single edge labelled by 
$1_x$ as shown below. 

\setlength{\unitlength}{1mm} 
\begin{center}
\begin{picture}(10,10)(2,0)
\put(5,0){\line(0,1){10}}    %
\put(5,10){\vector(0,-1){6}}  %
\put(6,5){\makebox(0,0)[l]{$1_x$}}  %
\end{picture}
\end{center} The source of $\alpha$ is empty, so we can define $g_0$ by
    \[g_0([x]) = [(\alpha_x,\cdot)]\]
where $(\alpha_x, \cdot) \in \elt{Q^+}\times_{\free{S'}} 
\free{A}$, and observe immediately that \[x \cong x' \in o(Q) \iff 
1_x \cong 1_{x'} \in \elt{Q}.\]  Furthermore we have 
    \[d_0[x] = \mu_T[x] = [1_x] = f_2 g_0 [x]\]
as required. 

For the induction step, suppose we have constructed $g_k$ 
satisfying the commuting condition; we seek to construct 
    \[g_{k+1} : C^{(k+1)} {\hookrightarrow} A_2\]
satisfying the condition. Consider 
    \[y \in C^{(k+1)} = S \ \amalg \ T(C^{(k)}) \times_{S^\prime} A.\]
If $y \in S$ then put $g_{k+1}(y) = g_0(y)$. Otherwise, we have 
\[y = (\alpha,a) \in T(C^{(k)}) \times_{S^\prime} A.\]  Here the map
$T(C^{(k)}) \lra S'$ is given by $Tf^{(k)}$.  Recall that by 
definition of $T$, $T(C^{(k)})$ is equivalent to the pullback 

\begin{center}
\begin{picture}(50,35)(5,5)        %

\put(10,10){\makebox(0,0){$\elt{Q}$}}  %  bottom left
\put(10,35){\makebox(0,0){$\cdot$}}  %  top left
\put(47,35){\makebox(0,0){$\free{(C^{(k)})}$}}  %  top right
\put(46,10){\makebox(0,0){$\free{S}$}}  %  bottom right

\put(13,35){\vector(1,0){25}}  %  top
\put(15,10){\vector(1,0){25}}  %  bottom
\put(10,32){\vector(0,-1){19}} %  left
\put(45,32){\vector(0,-1){19}} %  right

\put(8,23){\makebox(0,0)[r]{$$}} %  left
\put(47,23){\makebox(0,0)[l]{$\free{(f^{(k)})}$}}%  right
\put(27,37){\makebox(0,0)[b]{$$}}%  top
\put(27,8){\makebox(0,0)[t]{$$}} %  bottom

\end{picture}
\end{center} So, an element of \tck\ is an isomorphism class of arrows of $Q$
source-labelled by compatible elements of $C^{(k)}$. We write the 
pullback as \catck.  Then $Tf^{(k)}$ is the map given by the 
composite 
    \[\tck \map{\sim} \catck \lra \elt{Q} \map{\sim} S'.\]

Informally, $Tf^{(k)}$ removes the labels, leaving only the (\iso\ 
class of the) underlying arrow of $Q$. 

%9

Now we in fact exhibit a full and faithful functor 
    \[\catck \times_{S'} A \lra \xx{Q^+}{S'}{A}.\]
Let $((\beta,\underline{b}),a) \in \catck \times_{S'} A$.  So 
$\beta \in \elt{Q}$, $\underline{b} = \xv{b} \in \free{(C^{(k)})}$ 
and $a \in A$ such that $[s_Q(\beta)] = (f^{(k)}(b_1), \ldots, 
f^{(k)}(b_n))$ and $f(a) = [\beta]$. 

Informally, we have an arrow $\beta$ of Q, source-labelled by the 
$b_i \in \ck$, and a compatible label $a \in A$. We seek a formal 
composite of labelled arrows, of depth up to $k+1$.  By induction, 
we already have for each element of \ck\ a formal composite of 
labelled arrows, of depth up to $k$.  So we aim to form a formal 
composite of these together with $\beta$ labelled by $a$. 

By induction we have for each $1 \leq i \leq n$ 
    \[g_k(b_i) = (\pi_i,p_i) \in \xx{Q^+}{S'}{A}.\]
The commuting condition implies that for each $i$ 
    \[[s_Q(\beta)_i] = [t_Q t_{Q^+} (\pi_i)].\]
This gives us a way of constructing a new element of \elt{Q^+} 
from the data given, since each $\pi_i$ can be composed with 
$\beta$ at the $i$th place, via the appropriate 
object-isomorphism. That is, we form a tree by induction, as shown 
in the following diagram 

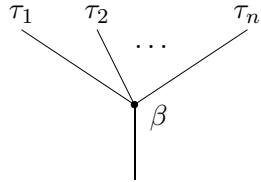
\setlength{\unitlength}{0.5mm} 
\begin{center}
\begin{picture}(80,55)      % tree with centre point

\put(40,0){\line(0,1){20}}   %
\put(40,20){\line(-3,2){30}} %
\put(40,20){\line(-1,2){10}} %
\put(40,20){\line(3,2){30}}  %
\put(40,20){\circle*{1.7}}     %

\put(10,43){\makebox[0pt]{$\tau_1$}}  % label 1
\put(30,43){\makebox[0pt]{$\tau_2$}}  % label 2
\put(70,43){\makebox[0pt]{$\tau_n$}}  % label n

\put(44,20){\makebox(0,0)[tl]{$\beta$}} % middle label

\put(40,35){$\ldots$} 
\end{picture}
\end{center} where $\tau_i$ is the tree for $\pi_i$. Each $\pi_i$
has its nodes (that is, source elements) labelled by elements of 
$A$; to complete the definition it remains only to `label' the 
node corresponding to $\beta$. But we have $f(a) = [\beta]$, that 
is, $a$ is a compatible label for $\beta$. So we let $a$ be the 
label for $\beta$. 

So we have defined a full and faithful functor 
    \[\catck \times_{S'} A \lra \xx{Q^+}{S'}{A}\]
inducing, on \iso\ classes, an embedding 
    \[g_{k+1}: \ck \hookrightarrow A_2\]
as required. We now check the commuting condition.  Informally, 
$d_k$ acts by ignoring the labels and composing the underlying 
arrows of $Q$, as does $\mu$. Since $\mu$ is induced from 
composition in $Q$, and $t_{Q^+}$ is constructed from composition 
of a formal composite of arrows of $Q$, we have $f_2 \circ g_{k+1} 
= d_{k+1}$ as required. 

So we have for each $k\ge 0$ an embedding $g_k$ as required. The 
$g_k$ then induce a map $A_1 \lra A_2$.  It is straightforward to 
check that this is surjective; by construction it makes the 
following diagram commute 

\setlength{\unitlength}{1mm} 
\begin{center}
\begin{picture}(40,25)(0,10)       %

\put(20,12){\makebox(0,0){$S'$}}  %  bottom
\put(5,30){\makebox(0,0){$A_1$}}  %  top left
\put(35,30){\makebox(0,0){$A_2$}}  %  top right

\put(9,30){\vector(1,0){22}}  %  top
\put(5,27){\vector(1,-1){13}}  %  left
\put(34,27){\vector(-1,-1){13}} %  right

\put(30,20){\makebox(0,0)[l]{$$}} %  right
\put(9,20){\makebox(0,0)[r]{$$}}%  left
\put(20,31){\makebox(0,0)[b]{$$}}%  top

\end{picture}
\end{center}so we have an isomorphism
    \[(A_1,f_1) \cong (A_2,f_2)\]
as required. 

Finally we check that the naturality condition for a monad 
opfunctor is satisfied.  Given a morphism $(A,f) \lra (B,g) \in 
\set/S'$ it is clear from the constructions that the following 
diagram commutes in $\cat{Set}/{S'}$ 

\setlength{\unitlength}{1mm} 
\begin{center}
\begin{picture}(40,28)(5,8)       %

\put(10,10){\makebox(0,0){$B_1$}}  %  bottom left
\put(10,30){\makebox(0,0){$A_1$}}  %  top left
\put(40,30){\makebox(0,0){$A_2$}}  %  top right
\put(40,10){\makebox(0,0){$B_2$}}  %  bottom right

\put(13,30){\vector(1,0){24}}  %  top
\put(13,10){\vector(1,0){24}}  %  bottom
\put(10,27){\vector(0,-1){14}} %  left
\put(40,27){\vector(0,-1){14}} %  right

\put(8,23){\makebox(0,0)[r]{$$}} %  left
\put(47,23){\makebox(0,0)[l]{$$}}%  right
\put(25,31){\makebox(0,0)[b]{$\simeq$}}%  top
\put(25,11){\makebox(0,0)[b]{$\simeq$}} %  bottom

\put(13,12){\vector(3,2){9}}  %    bl
\put(13,28){\vector(3,-2){9}}  %   tl
\put(37,28){\vector(-3,-2){9}}  %  tr
\put(37,12){\vector(-3,2){9}}  %   br
\put(25,20){\makebox(0,0)[c]{$S'$}} %  middle label

\end{picture}
\end{center}and the other axioms for a monad opfunctor are easily
checked. So we have \[(\cspl,\tpl) \cong (\cspr,\tpr)\] as 
required. 
\end{prf}

\addcontentsline{toc}{section}{References}
\bibliography{bib0209}

\nocite{bd2}
\nocite{hmp2} \nocite{hmp3} \nocite{hmp4}

\nocite{bae1}

\nocite{kel2}

\nocite{ks1}

\nocite{bkp1}

\end{document}